\documentclass[12pt,a4paper]{article}
\usepackage{t1enc}
\usepackage[latin1]{inputenc}
\usepackage[english]{babel}
\usepackage{graphicx}
\usepackage{mathptmx}

\font\bBB=msbm10
\font\bBBB=msbm10 at 18pt

\def\bBR{\mbox{\bBB R}}

\def\bBZ{\mbox{\bBB Z}}

\def\bBBZ{\mbox{\bBBB Z}}
\begin{document}
%\pagestyle{empty}
%\begin{centering}{\large Seeing the \Large{\bBBZ} $ _2$ in Three-Dimensional Rotations}
%\end{centering}
\title {Touching the $\bBBZ _2$ in Three-Dimensional Rotations}
\author{Vesna Stojanoska,\thanks{American University in Bulgaria, 2700 Blagoevgrad, Bulgaria, vstotka@gmail.com, 
vns021@aubg.bg
} \
 Orlin Stoytchev\thanks{American University in Bulgaria, 2700 Blagoevgrad, Bulgaria and Institute for Nuclear Research, 1784 Sofia, Bulgaria, ostoytchev@aubg.bg} }
\date{}
\maketitle
\abstract{A simple and self-contained proof is presented of the well-known fact that the
fundamental group of $SO(3)$ is $\bBZ _2$, using a relationship between closed paths in $SO(3)$ and braids.}
\section{Introduction}
\par\smallskip
The three-dimensional rotations form a Lie group, usually denoted
by $SO(3)$. A fascinating feature of our world  is that the
fundamental group of $SO(3)$ is $\bBZ _2$, or in short $\pi
_1(SO(3))\cong \bBZ _2$. This means that all closed paths in
$SO(3)$ starting and ending at the identity fall into two homotopy
classes --- those that are homotopic to the constant path and
those that are not. Composing two paths from the second class
yields a path from the first class. This topological property of
the space parametrizing rotations can also be stated as follows
--- a complete rotation of an object (corresponding to a closed path in $SO(3)$) may or
may not be continuously deformable to the trivial motion (i.e., no
motion at all); the composition of two motions that are not
deformable to the trivial one gives a motion which is.\par A
manifold with a fundamental group $\bBZ _2$ is a challenge to our
imagination --- it is easy to visualize spaces with fundamental
group $\bBZ$, (the punctured plane), or $\bBZ \times \bBZ \cdots
\times \bBZ$, (plane with several punctures), or even $\bBZ\oplus
\bBZ$ (torus), but there is no two-dimensional manifold embedded
in $\bBR^3$ whose fundamental group is $\bBZ _2$.\par The peculiar
structure of $SO(3)$ plays a crucial role in our physical world.
We know that there are two quite different types of elementary
particles, bosons and fermions. The quantum state of a boson is
described by a (possibly multi-component) wave function, which
remains unchanged when a full ($360^{\circ}$) rotation of the
coordinate system is performed, while the wave function of a
fermion gets multiplied by $-1$ under a complete rotation.
Somewhat loosely speaking, this comes from the fact that only the
modulus of the wave function has a direct physical meaning and
therefore the wave function need not transform under a true
representation of $SO(3)$ but just under a projective
representation, which is a true representation of its universal
covering group $SU(2)$ \cite{Wigner, Bargmann}.\par The standard
way of showing that $\pi _1(SO(3))\cong \bBZ _2$ is to prove that
$SO(3)$ is homeomorphic to $S^3/i$, where $i$ denotes the
identification of diametrically opposite points of $S^3$. Once the
latter is known, it is quite easy to see that paths starting from
one pole of $S^3$ and ending at the other pole will be closed
paths in $S^3/i$, which are not contractible, but the composition
of two such paths gives a contractible path. In order to find the
topological structure of $SO(3)$ one normally uses Lie group and
Lie algebra theory. Namely, it is shown that the group $SU(2)$ of
unitary $2\times 2$ matrices with determinant 1 is homeomorphic to
$S^3$ and that there is a 2--1 homomorphism $SU(2)\to SO(3)$,
which is a local isomorphism, and which sends diametrically
opposite points in $SU(2)$ to the same point in $SO(3)$. \par
There is a more direct way to exhibit the topological structure of
the parameter space of three-dimensional rotations. It only uses
the fact that rotations are represented by $3\times 3$ matrices
and any such matrix must have one (real) eigenvalue and a
corresponding eigenvector. Then one shows easily that this element
of $SO(3)$ must be a rotation around this eigenvector. In other
words any element of $SO(3)$ is a rotation around some axis and we
need to specify the angle of rotation and the orientation of that
axis in $\bBR^3$.\par In the present paper we present an
alternative way of proving that $\pi _1(SO(3))\cong \bBZ _2$. It
does not use Lie group theory or even matrices. It is purely
algebraic-topological in nature and very visual. It displays a
simple connection between full rotations (closed paths in $SO(3)$)
and braids. We believe that this may be an interesting way of
showing a nontrivial topological result to students in
introductory geometry and topology courses as well as a suitable
way of introducing braids and braid groups. \par The goal of this
paper is mostly pedagogical --- presenting in a self-contained and
accessible way a set of results that are basically known to
algebraic topologists and people studying braid groups. The fact
that the first homotopy group  of $SO(3)$ can be related to
spherical braids is a special case (in disguise) of the following
general statement \cite{Fad1}: ``The configuration space of three
points on an $r$-sphere is homotopically equivalent to the Stiefel
manifold of orthogonal two-frames in $r+1$-dimensional Euclidean
space''. Fadell \cite{Fad1} considers a particular element of $\pi
_1(SO(3))$ and uses the fact that it has order 2 to prove a
similar statement for a corresponding braid. Our direction is
opposite - we analyze braids to deduce topological properties of
$SO(3)$.\par In the next section we describe a simple experiment
that actually demonstrates the $\bBZ _2$ in three-dimensional
rotations. Then in section 3 we give a formal treatment of that
experiment. We construct a map from $\pi _1(SO(3))$ into a certain
factorgroup of a subgroup of the braid group with three strands.
We prove that this map is an isomorphism and that the image is
$\bBZ _2$.
\par
\section {The experiment}
\begin{figure} [h]
\setlength{\abovecaptionskip}{20pt}
\setlength{\belowcaptionskip}{0pt}
\centering
\includegraphics[scale=0.4
]{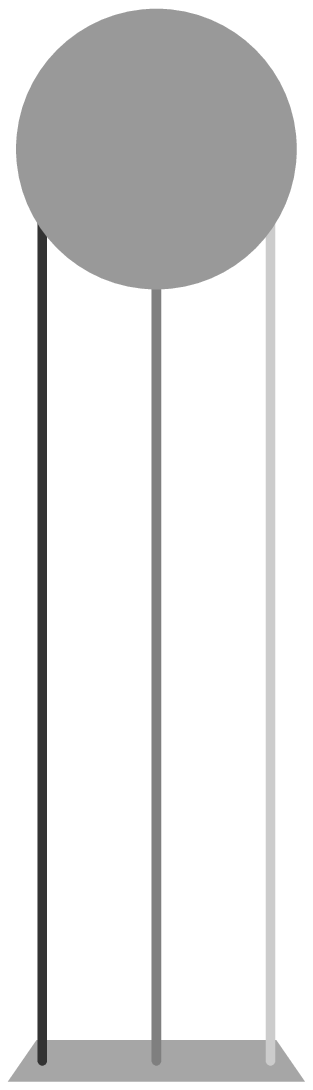}    \hfil         \includegraphics[scale=0.4]{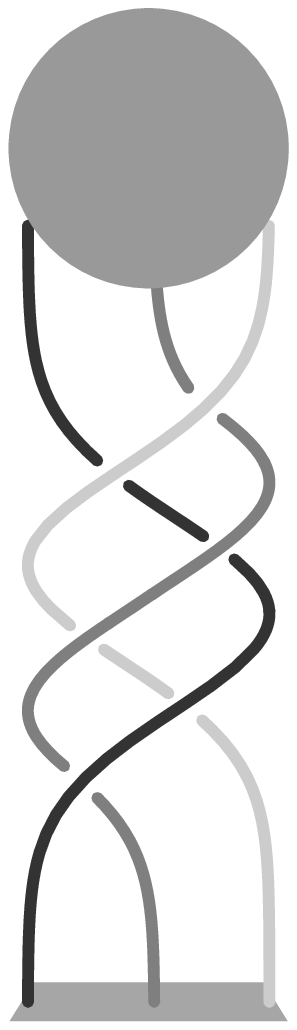}
\caption{Rotating a ball with strands attached}
\label{Fig1}
\end{figure}

Take a ball (tennis ball will do) and attach three strands to three different points on its surface. Attach the other ends of the strands to three different points on the surface of your desk  (Figure \ref{Fig1}).
\bigskip
Perform an arbitrary number of full rotations of the ball around arbitrary axes. You will get an entangled ``braid''. Now keep the orientation of the ball fixed. If the total number of full rotations is even, you can always untangle the ``braid'' by flipping strands around the ball. If the number of rotations is odd you will never be able to untangle it, but you can always reduce it to one simple configuration, e.g., the one obtained by rotating the ball around the first point and twisting the second and third strands around each other. \par
As is to be expected, rotations that can be continuously deformed to the trivial rotation (i.e., no rotation) lead to trivial braiding. At this point we can only conjecture from our experiment that the fundamental group of $SO(3)$ contains as a factor $\bBZ _2$. \par
\section{Relating three-dimensional rotations to braids}
\begin{figure} [h]
\setlength{\abovecaptionskip}{20pt}
\setlength{\belowcaptionskip}{0pt}
\centering
\includegraphics[scale=0.30]{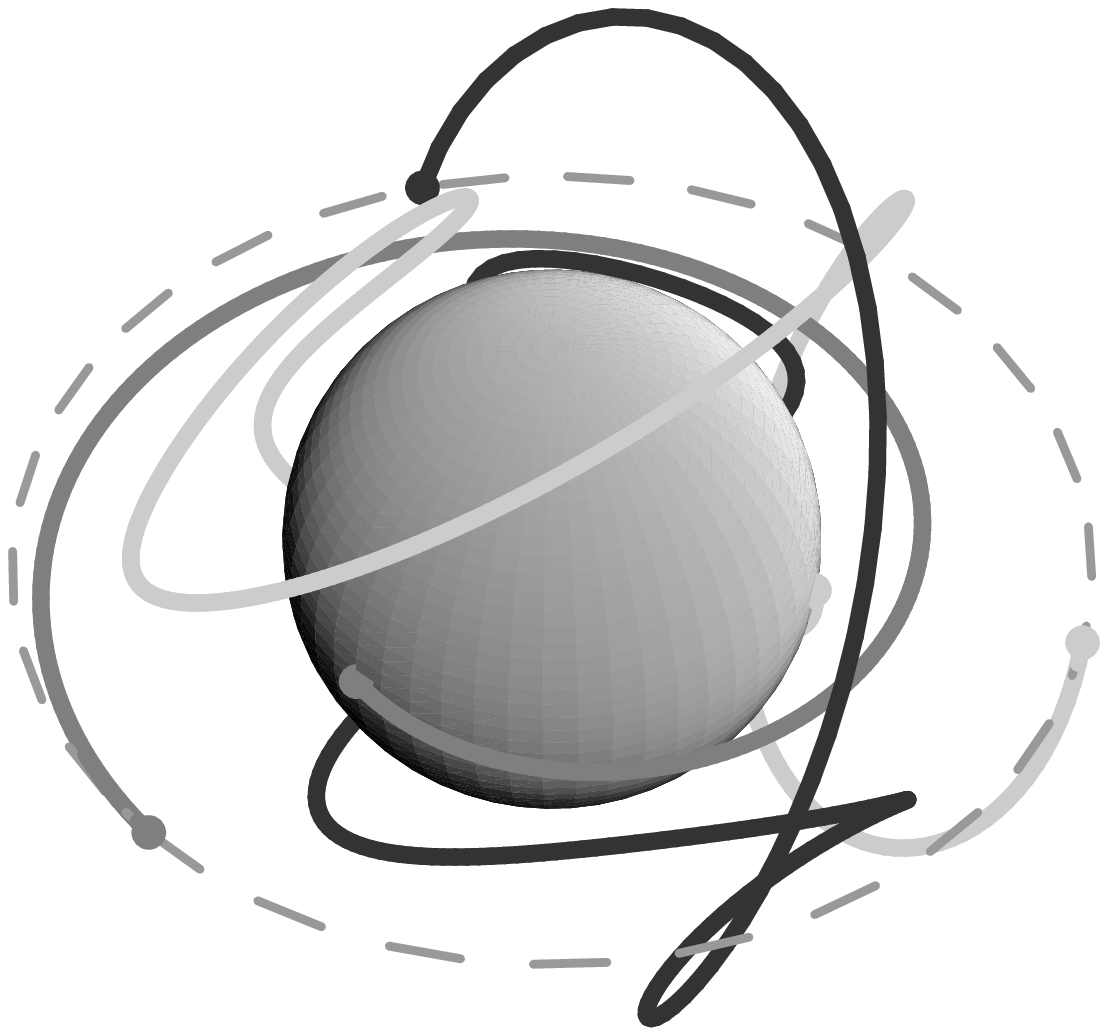}        \includegraphics[scale=0.30]{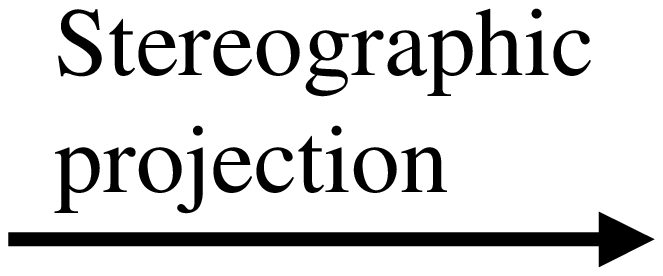}  \includegraphics[scale=0.50]{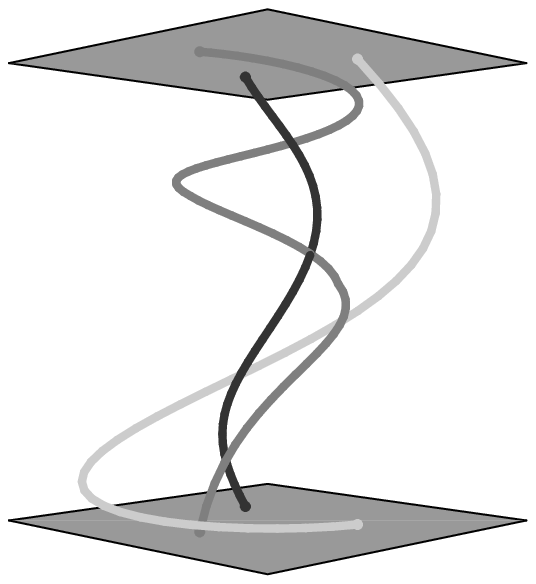}
\caption{A ``spherical braid'' and a normal braid}
\label{Fig2}
\end{figure}

With each closed path in $SO(3)$ we associate three closed paths in $\bBR^3$ starting at the sphere with radius 1 and ending at the sphere with radius 1/2. We may think of continuously rotating a sphere from time $t=0$ to time $t=1$ so that the sphere ends up with the same orientation as the initial one. Simultaneously we shrink the radius of the sphere from 1 to 1/2 (see Figure \ref{Fig2}).
\bigskip
Any three points on the sphere will trace three continuous paths in $\bBR^3$, which do not intersect each other. Furthermore, for fixed $t$ the three points on these paths lie on the sphere with radius $1-t/2$.
To formalize things, let $\omega(t), \ \ t\in[0,1]$ be any continuous path in $SO(3)$ with $\omega(0)=\omega(1)=I$.
$\omega(t)$ acts on vectors (points) in $\bBR^3$.
Take three initial points in $\bBR^3$, e.g., $\hbox{\bf x}_0^1=(1,0,0)$, $\hbox{\bf x}_0^2=(-1/2,\sqrt 3 /2,0)$,
$\hbox{\bf x}_0^3=(-1/2,-\sqrt 3 /2,0)$. Define three continuous paths by
$$
\hbox{\bf x}^i(t):=(1-t/2)\omega (t)(\hbox{\bf x}_0^i),\quad t\in [0,1],\quad i=1,2,3\ \ .
$$
 In this way we get an object that will be called a {\it spherical braid} --- several distinct points on a sphere and the same number of points, in the same positions, on a smaller sphere, connected by strands in such a way that the radial coordinate of each strand is monotonic in $t$.
\par\smallskip\noindent
{\bf Note} One can multiply two spherical braids by connecting the ends of the first to the beginnings of the second (and rescaling the parameter). When one considers classes of isotopic spherical braids one obtains the so called {\it braid group of the sphere} \cite{Fad2}, which algebraically is $B_3/R$ (see below). This is known as the mapping-class group of the sphere (with 0 punctures and 0 boundaries) and has been studied by topologists.
\par\smallskip
We can map our spherical braid to a conventional one using stereographic projection (Figure \ref{Fig2}). First we choose a ray starting at the origin and not intersecting any strand. Then we map stereographically any point on a sphere with radius $1/2 \leq \rho \leq 1$, except the point where the ray intersects that sphere and with respect to which we project, to a point in a corresponding horizontal plane. Finally we define the $z$-coordinate of the image to be $z=-\rho$ .\par
\begin{figure}
\setlength{\abovecaptionskip}{20pt}
\setlength{\belowcaptionskip}{20pt}
\centering
\includegraphics[scale=0.35]{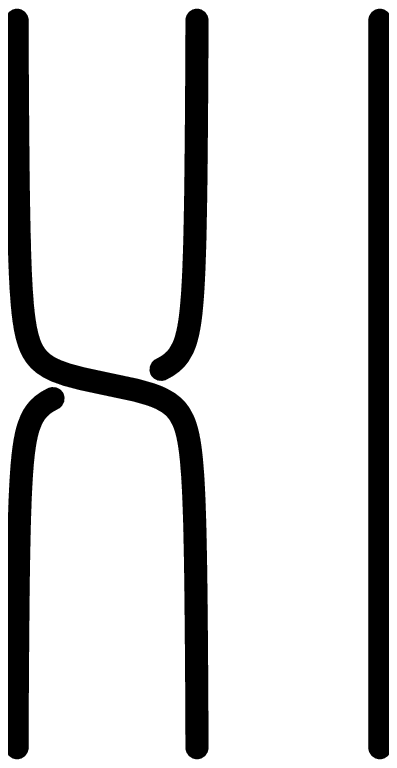}  \hskip80pt  \includegraphics[scale=0.35]{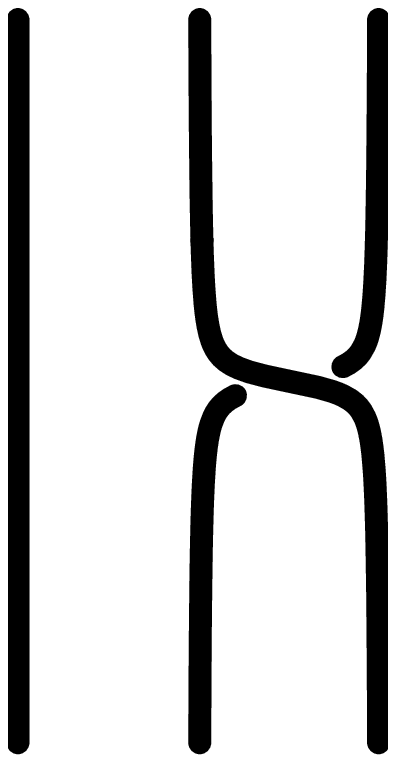}

\caption{The generators $\sigma_1$ and $\sigma_2$ of $B_3$ }
\label{Fig3}
\end{figure}

Recall the usual notion of {\it braids}, introduced by Artin \cite{Artin}. (See also \cite{Bir} for a contemporary review of the theory of braids and its relations to other subjects.) One takes two planes in $\bBR^3$, let's say parallel to the
$XY$ plane, fixes $n$ distinct points on each plane and connects each point on the lower plane with a point on the upper plane by a continuous path (strand). The strands do not  intersect each other. In addition the $z$-coordinate of each strand is a monotonic function of the parameter of the strand and thus $z$ can be used as a common parameter for all strands. Two different braids are considered equivalent or {\it isotopic} if there exists a homotopy of the strands (keeping the endpoints fixed), so that for each value of the homotopy parameter $s$ you get a braid, for $s=0$ you get the initial braid and for $s=1$ the final one. When the points on the lower and the upper plane have the same positions (their $x$ and $y$ coordinates are the same), we can multiply braids by stacking one on top of the other. Considering classes of isotopic braids with the multiplication just defined, the {\it braid group} is obtained. Artin showed that the braid group $B_n$ on $n$ strands has a presentation with $n-1$ generators and a simple set of relations - Artin's braid relations. We give them for the case $n=3$ since this is the one we are mostly interested in. In this case the braid group $B_3$ is generated by the generators $\sigma_1$, corresponding to twisting of the first and the second strands, and $\sigma_2$, corresponding to twisting of the second and the third strands (the one to the left always passing behind the one to the right) (Figure \ref{Fig3}).
\bigskip
These generators are subject to a single braid relation (Figure \ref{Fig4}):
\begin{figure} 
\setlength{\abovecaptionskip}{20pt}
\setlength{\belowcaptionskip}{0pt}
\centering
\includegraphics[scale=0.40]{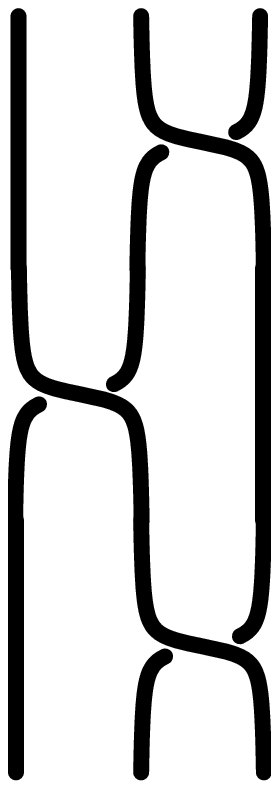} \quad\quad  \includegraphics[scale=0.20]{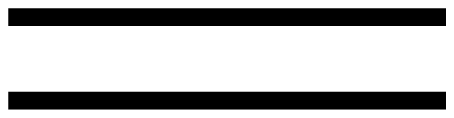}\quad\quad
   \includegraphics[scale=0.40]{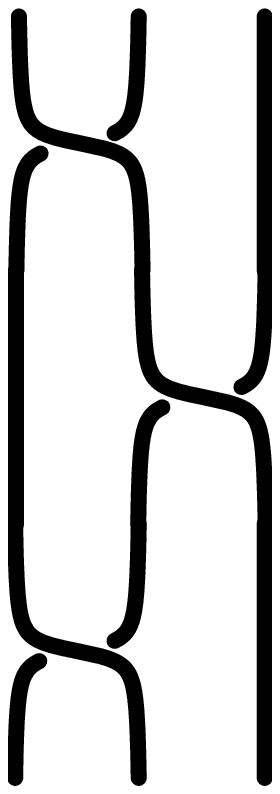}
\caption{The braid relation for $B_3$ }
\label{Fig4}
\end{figure}
\bigskip
\begin{equation}
\sigma_2\sigma_1\sigma_2=\sigma_1\sigma_2\sigma_1
\label{Artinrel}
\end{equation}
We say that $B_3$ has a {\it presentation} with generators $\sigma_1$ and $\sigma_2$ and defining relation given by Equation \ref{Artinrel}, or in short:
\begin{equation}
B_3=\langle \sigma_1,\,\sigma_2\,;\sigma_1\sigma_2\sigma_1\sigma_2^{-1}\sigma_1^{-1}\sigma_2^{-1}\rangle
\end{equation}
\par
\begin{figure} [h]
\setlength{\abovecaptionskip}{20pt}
\setlength{\belowcaptionskip}{0pt}
\centering
\includegraphics[scale=0.45]{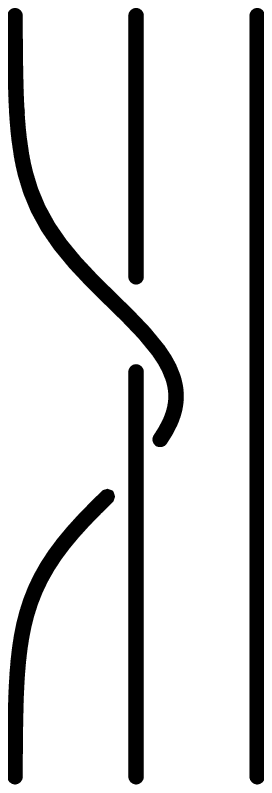}  \hskip60pt \includegraphics[scale=0.37]{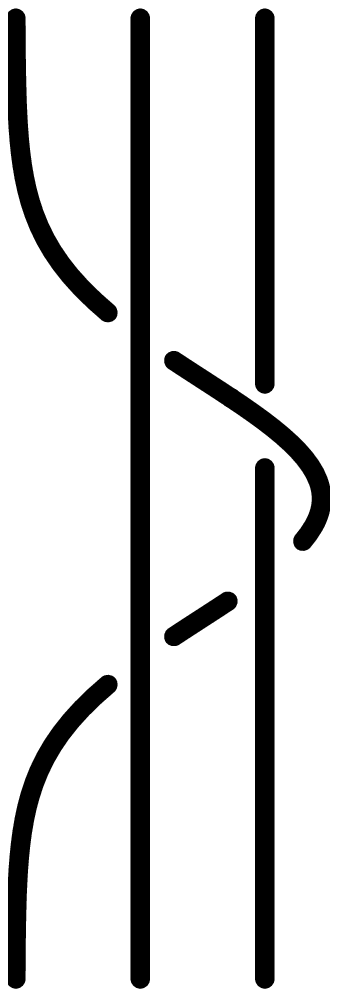}    \hskip60pt
\includegraphics[scale=0.4]{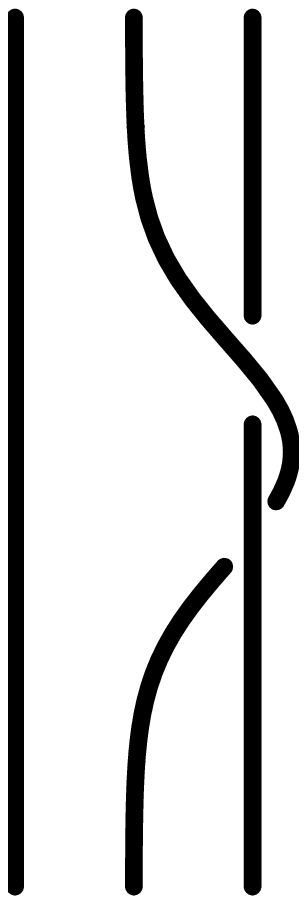}
\caption{The generators $a_{12}$, $a_{13}$ and $a_{23}$ of $P_3$ }
\label{Fig5}
\end{figure}

In our case, since a full rotation of the sphere returns the three
points to their original positions, we always get {\it pure
braids}, i.e., braids for which any strand connects a point on the
lower plane with its translate on the upper plane. Pure braids
form a subgroup of $B_3$ which is denoted by $P_3$. Note that
intuitively there is a homomorphism $\pi$ from  $B_3$ to the
symmetric group $S_3$ since any braid from $B_3$ permutes the
three points. Formally one defines $\pi$ on the generators by
\begin{equation}
\pi(\sigma_1)(1,2,3)=(2,1,3),\quad \pi(\sigma_2)(1,2,3)=(1,3,2)
\end{equation}
and then extends it to the whole group $B_3$ (it is important that $\pi$ maps Equation \ref{Artinrel} to the trivial identity). Pure braids are precisely those that do not permute the points and therefore we can give the following algebraic characterization of $P_3$:
$$P_3:=Ker\, \pi$$
Alternatively, $S_3$ is the quotient of $B_3$ by the additional
equivalence relations $\sigma_i^2\sim I, \ i=1,2$ and if $N$ is
the minimal normal subgroup containing $\sigma_i^2$, then
$\pi:B_3\to B_3/N$ is the natural projection. It is then easy to
see that the kernel of $\pi$ has to be a product of words of the
following type:
$$
\sigma_{i_1}^{\pm 1}\sigma_{i_2}^{\pm 1}\cdots\sigma_{i_k}^{\pm 1}\sigma_{i_{k+1}}^{\pm 2}\sigma_{i_k}^{\pm 1}
\cdots\sigma_{i_2}^{\pm 1}\sigma_{i_1}^{\pm 1}
$$
The whole subgroup $P_3$ can in fact be generated by the following three {\it twists} (Figure \ref{Fig5})
\begin{equation}
a_{12}:=\sigma_1^2,\quad a_{13}:=\sigma_2\sigma_1^2\sigma_2^{-1}=\sigma_1^{-1}\sigma_2^2\sigma_1,\quad a_{23}:=\sigma_2^{2}
\label{gener}
\end{equation}
\par
In our construction so far we mapped any closed path in $SO(3)$ to a spherical braid and then, using stereographic projection, to a conventional pure braid. The last map, however, depends on a choice of a ray in $\bBR^3$ and, what is worse, spherical braids that are isotopic (in the obvious sense) may map to nonisotopic braids. To mend this, we will identify certain classes of braids in
$P_3$. Namely, we introduce the following equivalence relations (see Figure \ref{Fig6}):
\begin{figure} [h]
\setlength{\abovecaptionskip}{20pt}
\setlength{\belowcaptionskip}{0pt}
\centering
\includegraphics[scale=0.50]{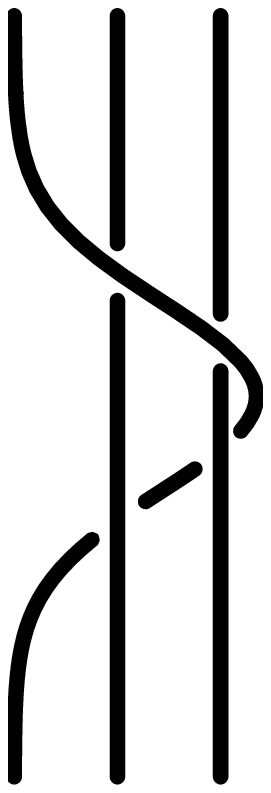}  \hskip40pt
\includegraphics[scale=0.45]{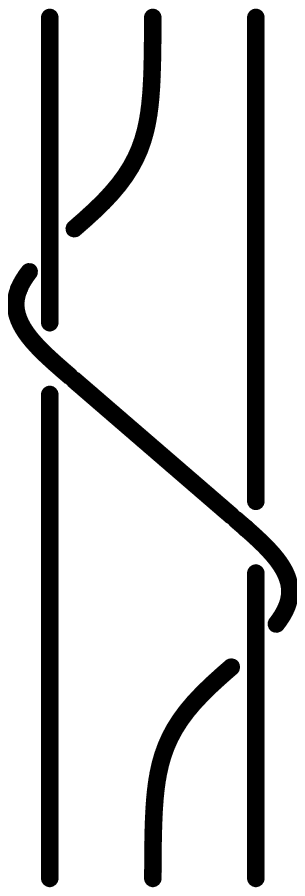}  \hskip40pt
\includegraphics[scale=0.50]{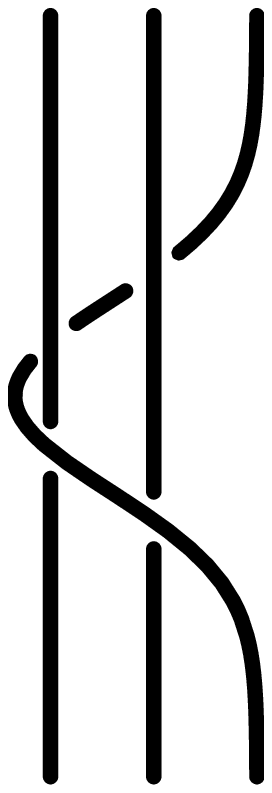}
\caption{The flips $r_1$, $r_2$ and $r_3$ }
\label{Fig6}
\end{figure}
\begin{equation}
r_1:=\sigma_1\sigma_2^2\sigma_1\sim I,\quad r_2:=\sigma_1^2\sigma_2^2\sim I,\quad
r_3:=\sigma_2\sigma_1^2\sigma_2\sim I. \label{r}
\end{equation}
On our model with the tennis ball from Section 2 $r_i,\,i=1,2,3$ correspond to {\it flips} of the $i$-th strand above and around the ball. Such motions  lead to isotopic spherical braids, as will be shown later.

\par\smallskip\noindent
{\bf Note} Strictly speaking, we have to make sure that our final result will not change when any strand crosses at any part of the braid the ray which we use for the stereographic projection. This means that we have to factorize by the normal closure in $B_3$ (not in $P_3$ !) of the generators $r_i,\ i=1,2,3$,  i.e., the smallest normal subgroup in $B_3$ containing these three generators. This would then allow us to set to $I$ any $r_i$ in any part of a word. One sees immediately that only one of the generators is needed then, since the other two will be contained in the normal closure of the first. We noticed, however, that we managed to untie any trivial braid just by a sequence of the three flips $r_i$ defined in Equation \ref{r} and their inverses, performed at the end of the braid. At the same time a nontrivial braid, corresponding to an odd number of rotations, cannot be untied even if we allow flips in any part of the braid. This can only be true if the flips $r_i$ generate a normal subgroup in $B_3$ (which of course then coincides with the normal closure of any of the $r_i$ and is also normal in $P_3$).
\par\bigskip\noindent
{\bf Lemma 1}\par\noindent
{\it The subgroup $R\subset P_3$, generated by $r_1$, $r_2$, $r_3$ is normal in $B_3$.}\par\smallskip\noindent
{\bf Proof:}
We need to show that we can represent all conjugates of $r_i$ with respect to the generators of $B_3$ and their inverses as products of the $r_i$ and their inverses. Straightforward calculations, using multiple times Artin's braid relation (Equation \ref{Artinrel}) give the following identities:
$$
\begin{array}{l ll}
\sigma_1r_1\sigma_1^{-1}=r_2, &\sigma_2r_1\sigma_2^{-1}=\sigma_2^{-1}r_1\sigma_2=r_1, \nonumber\\
\sigma_1r_2\sigma_1^{-1}=r_2r_1r_2^{-1}, &\sigma_2r_2\sigma_2^{-1}=r_3,\nonumber\\
\sigma_1r_3\sigma_1^{-1}=\sigma_1^{-1}r_3\sigma_1=r_3, &\sigma_2r_3\sigma_2^{-1}=r_1^{-1}r_2r_1=r_3r_2r_3^{-1},\nonumber\\
\sigma_1^{-1}r_1\sigma_1=r_1^{-1}r_2r_1, &\sigma_1^{-1}r_2\sigma_1=r_1, \nonumber\\
\sigma_2^{-1}r_2\sigma_2=r_1r_3r_1^{-1}=r_2^{-1}r_3r_2, \quad\quad &\sigma_2^{-1}r_3\sigma_2=r_2.\nonumber\\
\end{array}
$$
We demonstrate as an example the proof of the first identity in the second line. We have
\[\sigma_{1}\sigma_{2}\sigma_{1} = \sigma_{2}\sigma_{1}\sigma_{2} \]
\[\sigma_{2}\sigma_{1}\sigma_{2}\sigma_{1}\sigma_{2}\sigma_{1} = \sigma_{2}^{2}\sigma_{1}\sigma_{2}^{2}\sigma_{1} \]
\[\sigma_{1}\sigma_{2}\sigma_{1}\sigma_{2}\sigma_{1}\sigma_{2} = \sigma_{2}^{2}\sigma_{1}\sigma_{2}^{2}\sigma_{1} \]
\[\sigma_{1}\sigma_{2}^{2}\sigma_{1}\sigma_{2}^{2} = \sigma_{2}^{2}\sigma_{1}\sigma_{2}^{2}\sigma_{1} \]
\[\sigma_{1}\sigma_{2}^{2}\sigma_{1} = \sigma_{2}^{2}\sigma_{1}\sigma_{2}^{2}\sigma_{1}\sigma_{2}^{-2} \]
\[\sigma_{1}^{3}\sigma_{2}^{2}\sigma_1^{-1} = \sigma_1^2\sigma_{2}^{2}\sigma_{1}\sigma_{2}^{2}\sigma_{1}\sigma_{2}^{-2}\sigma_1^{-2} \]

and therefore

\[\sigma_1r_2\sigma_1^{-1}=
\sigma_{1}\cdot \sigma_{1}^{2} \sigma_{2}^{2}\cdot \sigma_{1}^{-1} = \sigma_{1}^{2}\sigma_{2}^{2} \cdot \sigma_{1}\sigma_{2}^{2}\sigma_{1} \cdot \sigma_{2}^{-2}\sigma_{1}^{-2}=r_2r_1r_2^{-1}\]

\par\bigskip
By suitable full rotations we obtain all generators of $P_3$. For
example, $a_{12}$ is obtained by rotating around the vector
$\hbox{\bf x}_0^3=(-1/2,-\sqrt 3 /2,0)$ and it twists the first
and the second strand. Furthermore, homotopies between closed
paths in $SO(3)$ correspond to isotopies of the spherical braids
and thus homotopic closed paths in $SO(3)$ will be mapped to the
same element in the factorgroup $P_3/R$.\par\bigskip\noindent {\bf
Proposition 1}\par\noindent {\it The factorgroup $P_3/R$ is
isomorphic to $\bBZ_2$.}\par\smallskip\noindent {\bf Proof:} To
make notations simpler we use the same letter to denote both a
representative of a class in $P_3/R$ and the class itself,  hoping
that the meaning is clear from the context. In $P_3/R$ we have
$$\sigma_1\sigma_2^2=\sigma_1^{-1}=\sigma_2^2\sigma_1,$$
 and
$$\sigma_2\sigma_1^2=\sigma_2^{-1}=\sigma_1^2\sigma_2.$$
Now we obtain the following sequence of identities:
$$
\begin{array}{l l}
\sigma_2\sigma_1^2=\sigma_1^2\sigma_2, &\sigma_1\sigma_2\sigma_1^2=\sigma_1^3\sigma_2, \nonumber \\
\sigma_2\sigma_1\sigma_2\sigma_1=\sigma_1^3\sigma_2, &\sigma_1\sigma_2\sigma_1\sigma_2\sigma_1=\sigma_1^4\sigma_2, \nonumber \\
\sigma_1\sigma_2^2\sigma_1\sigma_2=\sigma_1^4\sigma_2, &I=\sigma_1^4. \nonumber
\end{array}
$$
We have used twice the braid relation (Equation \ref{Artinrel}) and the first equivalence relation in Equation \ref{r}.
In completely the same way we prove
$$
\sigma_2^4=I.
$$
Combining the last two results with the equivalence relations (Equation \ref{r}) we finally get
\begin{equation}
\sigma_1^2=\sigma_1^{-2}=\sigma_2^2=\sigma_2^{-2}
\label{sigma}
\end{equation}
It is now clear that in $P_3/R$ the three generators, defined in Equation \ref{gener} reduce to one element of order 2. Therefore they generate $\bBZ_2$. This completes the proof.\par\smallskip
So far we have constructed a map $\pi_1(SO(3))\to P_3/R$, which is onto by construction, and we have shown that the image is isomorphic to $\bBZ_2$. We now show that this map is actually an isomorphism.
\par\bigskip\noindent
{\bf Proposition 2}\par\noindent {\it The map $\pi_1(SO(3))\to
P_3/R$ is a monomorphism.}\par\smallskip\noindent {\bf Proof:} We
have to show that if a closed continuous path in $SO(3)$ is mapped
to a braid in $R$, then this path is homotopic to the constant
path. The proof basically reduces to the following observation ---
any spherical braid which is pure (the strands connect each point
on the outer sphere with the same point on the inner sphere)
determines a closed path in $SO(3)$. Two isotopic spherical pure
braids determine homotopic closed paths in $SO(3)$. Indeed, recall
that for a spherical braid we can parametrize the points on each
strand with a single parameter $t$ and that for a fixed $t$ all
three points lie on a sphere with radius $1-t/2$. These three
ordered points $\hbox{\bf x}^i(t), i=1,2,3$ give for every fixed
$t$ a nondegenerate triangle, oriented somehow in $\bBR^3$. Let
$\hbox{\bf l}(t)$ be the vector, connecting the center of mass of
the triangle with the vertex $\hbox{\bf x}^1(t)$, i.e., $\hbox{\bf
l}(t)=\hbox{\bf x}^1-(\hbox{\bf x}^1(t)+\hbox{\bf
x}^2(t)+\hbox{\bf x}^3(t))/3$ and define $\hbox{\bf
e}^1(t):=\hbox{\bf l}(t)/||\hbox{\bf l}(t)||$. Let $\hbox{\bf
e}^3(t)$ be the unit vector, perpendicular to the plane of the
triangle, in a positive direction relative to the orientation
(1,2,3) of the boundary. Finally, let $\hbox{\bf e}^2(t)$ be the
unit vector, perpendicular to both $\hbox{\bf e}^1(t)$ and
$\hbox{\bf e}^3(t)$, so that the three form a right-handed frame.
Then there is a unique element $\omega (t)\in SO(3)$ sending the
vectors $\hbox{\bf e}_0^1=(1,0,0)$, $\hbox{\bf e}_0^2=(0,1,0)$,
$\hbox{\bf e}_0^3=(0,0,1)$ to the triple $\hbox{\bf e}^i(t)$.
According to our definitions, $\omega(0)=\omega(1)=I$ and we get a
continuous function $\omega: [0,1]\to SO(3)$, where continuity
should be understood relative to some natural topology on $SO(3)$,
e.g., the strong operator topology.\par
   Now, if we have two isotopic spherical braids, by definition there are continuous functions $\hbox{\bf x}^i(t,s), i=1,2,3$ such that
$\hbox{\bf x}^i(t,s)$ is a braid for any fixed $s\in[0,1]$,
$\hbox{\bf x}^i(0,s)=\hbox{\bf x}_0^i$, $\hbox{\bf
x}^i(1,s)=\hbox{\bf x}_0^i /2$, $\hbox{\bf x}^i(t,0)$ give the
initial braid and $\hbox{\bf x}^i(t,1)$ give the final braid. By
assigning an element $\omega(t,s)$ to any triple $\hbox{\bf
x}^i(t,s)$ as described we get a homotopy between two closed paths
in $SO(3)$.\par Let $\omega'(t)$ be a closed path in $SO(3)$ which
is mapped to a braid $b$ in the class $r_1\in R$.  We can
construct a spherical braid, whose image is isotopic to that
braid. Let $\hbox{\bf z}$ be the point on the unit sphere with
respect to which we perform the stereographic projection. This can
always be chosen to be the north pole or a point very close to the
north pole (in case a strand is actually crossing the axis passing
through the north pole). Note that the points $\hbox{\bf x}_0^i,
i=1,2,3$ are on the equator. Construct a simple closed path on the unit sphere
starting and ending at $\hbox{\bf x}_0^1$ and going around
$\hbox{\bf z}$ in a negative direction (without crossing the
equator except at the endpoints). Thus we have two continuous
functions $\varphi(t),\ \theta(t), t\in [0,1]$ --- the spherical
(angular) coordinates describing this path. Let $\hbox{\bf
x}^1(t)$ be the point in $\bBR^3$ whose spherical coordinates are
$\rho(t):=1-t/2,\ \varphi(t),\ \theta(t)$ and let $\hbox{\bf
x}^i(t):=(1-t/2)\hbox{\bf x}_0^i, i=2,3$. These three paths give
the required spherical braid. It is isotopic to the trivial braid,
coming from the constant path in $SO(3)$ and at the same time it
is isotopic to the preimage of $b$ under the stereographic
projection. In this way we see that $\omega'(t)$ must be homotopic to
the constant path. Obviously similar argument holds with $r_1$
replaced by $r_2$ and $r_3$ or the inverses. Since any element in
$R$ is a product of these generators and since products of
isotopic braids give isotopic braids, this completes the
proof.\par
\section{Further results and generalizations}
There is an obvious generalization of some of the results of the
previous sections to the case $n>3$. The minimal number of strands
that is needed to capture the nontrivial fundamental group of
$SO(3)$ is $n=3$. When $n>3$ any full rotation will give rise to a pure
spherical braid but the whole group of pure braids will not be
generated in this way. It is relatively easy to see that in this
way, after projecting stereographically, one will obtain a
subgroup of $P_n$, generated by a single {\it full twist} $d$ of all strands around an external
point and a set of $n$ flips $r_i$:
$$\begin{array}{l}
d:=(\sigma_{n-1}\cdots \sigma_2\sigma_1)^n,\\
r_1:=\sigma_1\sigma_2\cdots \sigma_{n-2}\sigma_{n-1}^2\sigma_{n-2}\cdots \sigma_1,\\
r_2:=\sigma_1^2\sigma_2\cdots\sigma_{n-2}\sigma_{n-1}^2\sigma_{n-2}\cdots \sigma_2,\\
$$r_i:=\sigma_{i-1}\cdots\sigma_2\sigma_1^2\sigma_2\cdots\sigma_{n-2}\sigma_{n-1}^2\sigma_{n-2}\cdots\sigma_i,
\quad i=2,3,\ldots n-1,\\
r_n:=\sigma_{n-1}\sigma_{n-2}\cdots\sigma_2\sigma_1^2\sigma_2\cdots\sigma_{n-2}\sigma_{n-1}.\\
\end{array}
$$

\begin{figure} [h]
\setlength{\abovecaptionskip}{20pt}
\setlength{\belowcaptionskip}{0pt}
\centering
\includegraphics[scale=0.28]{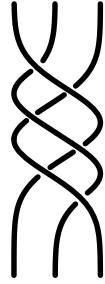}
\caption{The full twist $d$ in the case $n=3$ }
\label{Fig7}
\end{figure}

\begin{figure} [h]
\setlength{\abovecaptionskip}{20pt}
\setlength{\belowcaptionskip}{0pt}
\centering
\includegraphics[scale=0.45]{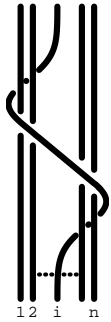}
\caption{The flip $r_i$ }
\label{Fig8}
\end{figure}

Figure
\ref{Fig7} shows a full twist for the case with 3 strands while Figure \ref{Fig8} shows a generic flip.
Straightforward calculations give the following generalization of Lemma 1:
\par\bigskip\noindent
{\bf Lemma 1$'$}\par\noindent
{\it The subgroup $R\subset P_n$, generated by $r_i$, $i=1,\ldots n$,  is normal in $B_n$.}\par\smallskip\noindent
{\bf Proof:} As in the proof of Lemma 1 we exhibit explicit formulas for the conjugates of all flips $r_i$:
$$
\begin{array}{l }
\sigma_j r_i \sigma_j^{-1}=\sigma_j^{-1}r_i \sigma_j=r_i,\quad i-j>1\ \hbox{or}\ j-i>0,\\
\sigma_{i-1}r_i\sigma_{i-1}^{-1}=r_i r_{i-1} r_i^{-1},\\
\sigma_{i-1}^{-1} r_i \sigma_{i-1}=r_{i-1},\\
\sigma_i r_i \sigma_i^{-1}=r_{i+1},\quad i<n-1\\
\sigma_i^{-1} r_i \sigma_i = r_i^{-1}r_{i+1} r_i,\quad i<n-1\\
\end{array}
$$
\par\bigskip
Let us denote by $S$ the subgroup, generated by $d$ and $r_i$.
Using purely topological information, namely that $\pi_1(SO(3))\cong\bBZ_2$, we can deduce the following
generalization of Proposition 1:
\par\bigskip\noindent
{\bf Proposition 1$'$}\par\noindent {\it The factorgroup $S/R$ is
isomorphic to $ \bBZ_2$.}\par\smallskip\noindent An equivalent
statement is that $d^2 \in R$.
\par It is tempting to try
generalizing the main result of this paper to higher-dimensional
rotations. A straightforward generalization fails since one would
produce braids in higher than three-dimensional space and these
can always be untangled. Note that in three dimensions we are able
to attribute a path in $SO(3)$ to any spherical braid with 3
strands but this is not the case for $n>3$ (4 points on $S^3$ may
not determine an orientation of the orthonormal frame in
$\bBR^4$.) It remains to be seen if a refinement of the methods
used can render meaningful results. 

\vfill\eject
\end{document}